\documentclass{article}
    \usepackage{amssymb, latexsym, amscd, amsthm, amstext, amsmath}

  \usepackage{amsopn}

  \newtheorem{thm}{Theorem}%[section]
 \newtheorem{lem}{Lemma}%[section]
 \newtheorem{cor}{Corollary}%[section]
 \newtheorem{rem}{Remark}%[section]
\begin{document}
 % Title. If the supplement option is on, then "Supplementary Material"
 % is automatically inserted before the title.
 \title{\textbf{Simultaneous diagonalization via congruence of $\lowercase{m}$ real symmetric matrices and its implications in optimization} \thanks{
 		Thi Ngan Nguyen  and Van-Bong Nguyen were funded by Tay Nguyen University under Grant  T2022-42CBTD.
 			Thi Ngan Nguyen was  funded by Vingroup Joint Stock Company and supported by the Domestic Master / PhD
 			Scholarship Programme of Vingroup Innovation Foundation (VINIF), Vingroup Big Data Institute
 			(VINBIGDATA), code VINIF. 2020.TS.75.}}
 
 % Authors: full names plus addresses.
 \author{Thi Ngan Nguyen\thanks{Department of Mathematics, Tay Nguyen
 		University, 632090, Vietnam
 		(\texttt{nguyenthingan@ttn.edu.vn, nguyenthingan@qnu.edu.vn}).}
 	\and Van-Bong Nguyen\thanks{Department of Mathematics, Tay Nguyen
 		University, 632090, Vietnam
 		(\texttt{nvbong@ttn.edu.vn}).}
 	\and Thanh-Hieu Le\thanks{Department of Mathematics and Statistics, Quy Nhon University, Vietnam
 		(\texttt{lethanhhieu@qnu.edu.vn}).}
 	\and Ruey-Lin Sheu\thanks{Department of Mathematics, National Cheng Kung
 		University, R.O.C. 70101, Taiwan
 		(\texttt{rsheu@mail.ncku.edu.tw}).}}

%\input{ex_shared}

 %Optional PDF information
%\ifpdf
%\hypersetup{
%  pdftitle={Simultaneous Diagonalization via Congruence of $m$ Real Symmetric Matrices},
%  pdfauthor={T. N. Nguyen, V. B. Nguyen, T. H. Le, and R. L. Sheu}
%}
%\fi

\maketitle

% REQUIRED
\begin{abstract}
Let $\{C_1, C_2, \ldots, C_m\},~m\ge2$ be a collection of $n\times n$ real symmetric matrices.
The objective of the paper is to offer an algorithm that finds a common congruence
matrix $R$ such that $R^TC_iR$ is real diagonal for every $C_i;$ or reports none of such kind.
The problem, referred to as the simultaneously diagonalization via congruence (SDC in short),
seems to be of pure linear algebra at first glance. However, for quadratically constrained quadratic programming (QCQP), if the quadratic forms are SDC, their joint range set is a closed convex polyhedral cone, which opens the possibility to extend the classical $\mathcal{S}$-lemma for more than two symmetric matrices.
In addition, under the SDC assumption of quadratic forms, QCQP can be recast in separable forms which is usually easier to tackle. It
is thus important to have a standard procedure for determining whether or not the SDC property holds for the underlined quadratic optimization problem. Our result solves a long standing problem posed by Hiriart-Urruty in 2007.
\end{abstract}

% REQUIRED
%\begin{keywords}

% REQUIRED
%\begin{MSCcodes}
%\end{MSCcodes}

\section{Introduction}
 Let $\mathcal{S}^n$ be the space of $n\times n$ real symmetric matrices and a collection $\mathcal{C}=\{C_1, C_2, \ldots, C_m\}\subset\mathcal{S}^n,~m\ge2$ be given. The problem to determine whether or not $\mathcal{C}$ is
SDC was first posed by Hiriart-Urruty and Torki \cite{Hiriart} in 2002 and later by Hiriart-Urruty \cite{Urruty} in 2007. The SDC problem remains open until now. The most current result was for $m=2$ by Jiang and Duan \cite{Jiang} in 2016. It was conjectured in \cite{Hiriart} and \cite{Urruty} that the SDC property can be used to ``\emph{generalize, to some extent, the Finsler-Calabi-type results}.'' In order to let the reader see what the
Finsler-Calabi-type results are all about and why they matter the optimization community, we feel that it is necessary to provide a short review.

Recall that classical Finsler theorem \cite{Finsler} in 1937 asserts that, for two symmetric matrices $C_1,C_2\in \mathcal{S}^n$ and $n\ge3$, the following two
statements are equivalent:
\begin{itemize}
\item[]$({\rm F_1})$~$(x^TC_1x=0) \wedge (x^TC_2x = 0)~\Longrightarrow~ (x=0);$
\item[]$({\rm F_2})$~($\exists \mu_1,\mu_2 \in \mathbb{R}$)~$\mu_1 C_1 + \mu_2 C_2\succ0.$
\end{itemize}
The result was later rediscovered by Calabi \cite{Calabi} in 1963.
Notice that both statements $({\rm F_1})$ and $({\rm F_2})$ imply, independently, that the collection $\mathcal{C}=\{C_1, C_2\}$ is SDC. Such an equivalence fails when the collection $\mathcal{C}$ consists of three or more symmetric matrices. In an attempt to search for a similar Finsler-Calabi-type result for $m\ge3,$ the SDC property of $\mathcal{C}$ does hold some hope.

Finsler-Calabi theorem, perhaps, was the earliest result that has the so-called ``hidden convexity,'' with which
certain non-convex quadratic optimization problems, under Slater's condition (or other types of constraint qualification), can be proved to adopt strong duality and solved from the dual side (usually by an SDP). Although $C_1,C_2\in\mathcal{S}^n$ may be indefinite, Brickman \cite{Brickman} in 1961 used a very enlightening
proof to show that, when $n\ge3,$ the joint range set
$$
\mathcal{R}(C_1,C_2)_{\|x\|=1}=\big\{ \big( x^TC_1x , x^TC_2x \big) \in \mathbb{R}^2 ~\big|~ \|x\|=1,~x \in \mathbb{R}^n \big\}
$$
is convex. Then, the non-trivial direction of Finsler-Calabi theorem ``$({\rm F_1})\Longrightarrow({\rm F_2})$'' has
become just to separate two convex sets:
$\mathcal{R}(C_1,C_2)_{\|x\|=1}$ and $\{0\}.$
See Section 2 in \cite{Hiriart} for detail. Moreover, the cone generated by the convex set
$\mathcal{R}(C_1,C_2)_{\|x\|=1}$ remains convex. That is,
\begin{equation}\label{Brickman}
Cone(\mathcal{R}(C_1,C_2)_{\|x\|=1})=\{\alpha\cdot  z ~|~\alpha\ge0, z\in \mathcal{R}(C_1,C_2)_{\|x\|=1}\}
\mbox{ is convex.}
\end{equation} However,
the convex cone $Cone(\mathcal{R}(C_1,C_2)_{\|x\|=1})$ is exactly the
joint range set of $(x^TC_1x , x^TC_2x)$ over the entire $\mathbb{R}^n.$ From \eqref{Brickman}, we immediately recover Dines's theorem \cite{Dines} in 1941 (for $n\ge3$) that
$$Cone(\mathcal{R}(C_1,C_2)_{\|x\|=1})=\mathcal{R}(C_1,C_2)_{\mathbb{R}^n}=\big\{\big( x^TC_1x , x^TC_2x \big) \in \mathbb{R}^2 \big|~x \in \mathbb{R}^n \big\} \mbox{ is convex.}$$
Dines' theorem later leads to the celebrated $\mathcal{S}$-lemma with which a bunch of fruitful results in quadratic optimization, especially the one-quadratic-constraint case, have been derived. Since those results are rather recent,
we only provide a few references for the interested reader. Please see 
 \cite[App. B]{Boyd},\cite{Polik-Terlaky07,JEYAKUMAR,Xia-Wang-Sheu16,YongXia,Bong}.

Brickman's result, from our understanding, is the most fundamental here. Both Finsler's theorem and Dines' theorem are its direct consequence. Furthermore, if we apply
Dines' theorem to Finsler's theorem, the latter becomes a type of $\mathcal{S}$-lemma. Specifically, the
following two statements are equivalent for $n\ge3.$
\begin{itemize}
\item[]$({\rm F'_1})$~$ x^TC_2x=0,~x\not=0~\Longrightarrow~ x^TC_1x>0;$
\item[]$({\rm F'_2})$~($\exists \mu\in\mathbb{R}$)~$C_1 + \mu C_2\succ0.$
\end{itemize}

The success of Finsler-Calabi theorem and the $\mathcal{S}$-lemma does not extend to more than two symmetric matrices. The major reason is that we no longer have the convexity of the joint range set
$$\mathcal{R}(C_1,C_2,C_3,\ldots,C_m)=\big( x^TC_1x , x^TC_2x, x^TC_3x, \ldots,x^TC_mx \big)\subset\mathbb{R}^m, (m\ge 3)$$
no matter the variable $x$ runs over the unit sphere or over $\mathbb{R}^n.$ However, under the assumption that the collection $\mathcal{C}$ is SDC, the aforementioned joint range set
$$\mathcal{R}(C_1,C_2,C_3,\ldots,C_m)_{\mathbb{R}^n}
=\big\{\big( x^TC_1x , x^TC_2x, x^TC_3x, \ldots,x^TC_mx \big)|x\in\mathbb{R}^n\big\}$$
is a closed convex polyhedral cone. In this case, the difficult NP-hard problem
\begin{equation*}
	{\rm(QCQP_H)} \hspace*{0.6cm}
	\begin{array}{lll}
	\lambda^*=&\min &x^TQ_1x\\
	&{\rm s.t.} & x^TQ_ix\le(=)~ b_i,~ i=2,\ldots,m,
	\end{array}
\end{equation*}
where $b_i\in\mathbb{R}$, is nicely reduced to a linear programming. In the final section of the paper after
the proof of our algorithm for determining the SDC of the collection $\mathcal{C}$ is clearly done, we will have a chance to loop back the issue for more discussions.

The type of problem ${\rm(QCQP_H)}$ has various applications, especially in the signal processing area.
Please refer to \cite{Luo} by Luo et al.. If linear terms are imposed, ${\rm(QCQP_H)}$ becomes the general (QCQP), in which many results set, as an assumption, the SDC of their underlined symmetric matrices.
For example,
Ben-Tal and Hertog \cite{Ben-Tal} showed that, if the two matrices in the objective and the constraint functions are SDC, a QCQP with one constraint (often known by the generalized trust region subproblem (GTRS)) can be equivalently transformed to a convex second-order cone problem (SOCP).
Jiang and Duan \cite{Jiang} apply the SDC of two matrices to solve the GTRS and its variants. Salahi and
Taati \cite{Salahi} also derived an efficient algorithm for solving GTRS under the SDC condition. Adachi and Nakatsukasa \cite{Adachi} use the SDC of two matrices to compute the positive definite interval of the matrix pencil and  propose an $O(n^3)$ novel eigenvalue-based algorithm for a definite feasible GTRS. Like the SDC
condition reduces ${\rm(QCQP_H)}$ to a linear programming, it also brings an equivalent and easier-to-tackle reformulation for the above applications; and provides a chance to solve the original problem in a much efficient way. This is certainly an indispensable benefit for large scale problems.

In fact, the SDC problem is very classical and has existed for a long time in the area of linear algebra. However, the progress seems to have stalled
at the case $m=2.$ Necessary and sufficient conditions for the SDC of two symmetric matrices $C_1,C_2$ can be found in \cite{Beck80} and also the references therein, but none is polynomially checkable. It has to wait
until Jiang and Li \cite{Jiang} in which an algorithm is developed to construct a congruence $R.$ Nevertheless, we find that the result of Jiang and Li \cite{Jiang} is not complete. A missing case not considered in their paper is now added to make it up in ours. Besides, we are able to extend their approach to answer the SDC of $\mathcal{C}=\{C_1, C_2, \ldots, C_m\}$ for any $m>2.$

We divide $\mathcal{C}=\{C_1, C_2, \ldots, C_m\}\subset\mathcal{S}^n$ into two cases. The first subcase is called the {\em nonsingular collection} (in Section 2), when at least one $C_i\in \mathcal{C}$ is non-singular. The other subcase is called the {\em singular collection} (in Section 3), when all $C_i's$ in $\mathcal{C}$ are non-zero but singular. When $\mathcal{C}$ is a nonsingular collection, we always assume that $C_1$ is non-singular. A non-singular collection will be denoted by $\mathcal{C}_{ns},$ while $\mathcal{C}_{s}$ represents the singular collection. Our main results in both singular as well as non-singular cases are iterative-based,
meaning that the arguments first apply to $\{C_1,C_2\},$ then apply for the second time on
$\{C_1,C_2,C_3\};$ the third time on $\{C_1,C_2,C_3,C_4\},$ and so forth.
However, in the proof, we only show the first few steps for people to see the recursive nature, rather than to attempt a complete set of mathematical induction, lest the messy notations may blur the essential concept of our approach. In the final section (Section 4), applications of the SDC for QCQP and for maximization of the sum of generalized Rayleigh quotients
and the extension of Finsler-Calabi type result for $m\ge3$ is discussed to finish the paper.

\section{The SDC problem of nonsingular collection}\label{sec01}

Consider a nonsingular collection $\mathcal{C}_{ns}=\{C_1, C_2, \ldots, C_m\}\subset \mathcal S^n$ and assume that $C_1$ is nonsingular. Let us outline the approach to determine the SDC of $\mathcal{C}_{ns}$. First, by
(i) and (iii) of Lemma \ref{L3.0} below, we show that if $\mathcal{C}_{ns}$ is SDC, it is necessary that
\begin{itemize}
  \item[(N1)] $C_1^{-1}C_i,~i=2,3,\ldots,m$ is real similarly diagonalizable;
  \item[(N2)] $C_j C_1^{-1} C_i$ is symmetric, for every $i = 2, 3, \ldots,m$ and every $j\not=i.$
\end{itemize}
Conversely, for the sufficiency, we use (N1) and (N2) to decompose, iteratively, all matrices in $\mathcal{C}_{ns}$ into block diagonal forms of smaller and smaller size until all of them become the so-called
non-homogeneous dilation of the same block structure (to be seen later) with certain scaling factors. Then, the SDC of $\mathcal{C}_{ns}$ is readily achieved.

The following technical lemmas from Linear Algebra play crucial roles.

\begin{lem}[Uhlig \cite{Uhlig76}]\label{bd2}
Let $K=\texttt{diag}(C(\lambda_{1}), \cdots, C(\lambda_{k}))$ be a Jordan matrix with eigenvalue $\lambda_{1},\cdots,\lambda_{k},$ where
$C(\lambda_{i})=\texttt{diag}(K_{i_{1}}(\lambda_{i}), K_{i_{2}}(\lambda_{i}),\cdots,K_{i_{t_i}}(\lambda_{i}))$ are Jordan blocks consisting of form
\[
K_{i_j}(\lambda_i)=
\left(
\begin{array}{ccccc}
 \lambda_{i} & 1 &  & & \\
0 & \lambda_{i} & 1& & \\
  & &\ddots&\ddots & \\
  &  & &\lambda_{i}&1\\
  &  & & & \lambda_{i}
\end{array}
\right)_{(i_j)}, \hspace{1cm}  j=1, 2, \cdots, t_i.
\]
Suppose, for a symmetric matrix $S,$  $SK$ is symmetric. Then, $S=\texttt{diag}(S_1,S_2,\ldots,S_k)$ is block diagonal with
${\rm dim} S_{i}={\rm dim }C(\lambda_{i}).$
\end{lem}

\begin{lem}\label{L3.0}
Let $C_1, C_2, C_3$ be real symmetric with $C_1$ non-singular. Suppose $P$ is non-singular such that $\hat{C}_1=P^TC_1P;~\hat{C}_2=P^TC_2P;~\hat{C}_3=P^TC_3P.$
%Then,
\begin{itemize}
  \item[(i)] If $C_1$ and $C_2$ are SDC, then $C_1^{-1}C_2$ is similarly diagonalizable.
  \item[(ii)] If $C_1^{-1}C_2$ is similarly diagonalizable, so is $\hat{C}_1^{-1}\hat{C}_2$.
  \item[(iii)] If $C_1,C_2,C_3$ are SDC, then $C_3C_1^{-1}C_2$ is symmetric.
  \item[(iv)] If $C_3C_1^{-1}C_2$ is symmetric, so is $\hat{C}_3\hat{C}_1^{-1}\hat{C}_2.$
\end{itemize}
\end{lem}
\begin{proof}
Since $C_1$ is non-singular, $\hat{C}_1$ is non-singular too. Moreover, there is
\begin{eqnarray}
  \hat{C}_1^{-1}\hat{C}_2 &=& P^{-1}C_1^{-1}C_2P; \label{easy-1}\\
  \hat{C}_3\hat{C}_1^{-1}\hat{C}_2 &=& P^TC_3C_1^{-1}C_2P.\label{easy-2}
\end{eqnarray}
For (i), suppose $C_1$ and $C_2$ are SDC, there is an invertible $Q$ such that
$$Q^TC_1Q=D_1,~Q^TC_2Q=D_2$$
where $D_1,~D_2$ are diagonal with $D_1$ invertible. By \eqref{easy-1},
$D_1^{-1}D_2=Q^{-1}C_1^{-1}C_2Q.$ Then, $C_1^{-1}C_2$ is similarly diagonalizable.

\noindent For (ii), suppose $C_1^{-1}C_2$ is similarly diagonalizable by a non-singular $H$, we can see that $P^{-1}C_1^{-1}C_2P$ is similarly diagonalizable by $P^{-1}H.$
By \eqref{easy-1}, $\hat{C}_1^{-1}\hat{C}_2$ is similarly diagonalizable.

\noindent For (iii), suppose $C_1,C_2,C_3$ are SDC. There is an invertible $Q$ such that
$$Q^TC_1Q=D_1,~Q^TC_2Q=D_2,~Q^TC_3Q=D_3$$
where $D_1,D_2,D_3$ are diagonal. By \eqref{easy-2}, $D_3D_1^{-1}D_2=Q^TC_3C_1^{-1}C_2Q,$ which shows that $C_3C_1^{-1}C_2$ is symmetric.

\noindent For (iv), suppose $C_3C_1^{-1}C_2$ is symmetric. By \eqref{easy-2},  $\hat{C}_3\hat{C}_1^{-1}\hat{C}_2$ is symmetric, too.
\end{proof}

In general, from (i) and (iii) of Lemma \ref{L3.0}, we find that, for a non-singular collection $\mathcal{C}_{ns}$ to be SDC, it is necessary that
\begin{itemize}
  \item[(N1)] $C_1^{-1}C_i,~i=2,3,\ldots,m$ is real similarly diagonalizable;
  \item[(N2)] $C_j C_1^{-1} C_i$ is symmetric, for every $i = 2, 3, \ldots,m$ and every $j\not=i.$
\end{itemize}
By Theorem \ref{dl80} and Theorem \ref{thm1} below, we will show that (N1) and (N2) are indeed sufficient for $\mathcal{C}_{ns}$ to be SDC. Let us begin with Lemma \ref{bd27}.

\begin{lem}\label{bd27}
Let $\mathcal{C}_{ns}=\{C_1, C_2, \ldots, C_m\}\subset \mathcal S^n$ be a non-singular collection with $C_1$ invertible. Suppose $C_1^{-1}C_2$ is real similarly diagonalized by invertible matrix $Q$ to have $r$ distinct eigenvalues $\beta_1,\ldots,\beta_r;$ each of multiplicity $m_t,~t=1,2,\ldots,r,$ respectively. Then,
\begin{eqnarray}
Q^TC_1Q &=& \texttt{diag}\underbrace{((A_1)_{m_1}, (A_2)_{m_2}, \ldots, (A_r)_{m_r})}_{m_1+\cdots+m_r=n,\small{\hbox{ each }A_t:\hbox{ sym. invert.}}}; \label{ct1}\\
Q^TC_2Q &=& \texttt{diag} (\beta_1A_1, \beta_2A_2 \ldots,\beta_rA_r).\label{ct1.1}
\end{eqnarray}

In addition, if $C_jC_1^{-1}C_2, j= 3, 4,\ldots, m,$ are symmetric, we can further block diagonalize $C_3, C_4, \ldots, C_m$ to adopt the same block structure as in \eqref{ct1}, such that
\begin{equation}
Q^TC_jQ=\texttt{diag} \underbrace{((C_{j1})_{m_1}, (C_{j2})_{m_2},\ldots, (C_{jr})_{m_r})}_{\small{\hbox{ each }}C_{jt}:\small{\hbox{ sym. }}}  \enskip j=3,4,\ldots,m.\label{ct1.2}
\end{equation}
\end{lem}

\begin{proof}
Since $C_1^{-1}C_2$ is similarly diagonalizable by $Q$, by assumption, there is
\begin{equation}\label{Lemma2-0}
J:= Q^{-1} C_1^{-1}C_2 Q =\texttt{diag}(\beta_1I_{m_1}, \ldots, \beta_r I_{m_r})
\end{equation}
with $m_1+m_2+\cdots+m_r=n.$ From \eqref{Lemma2-0}, we have, for $j= 1,2,\ldots, m,$
\begin{equation}\label{Lemma2-2}
(Q^TC_j Q)J= (Q^TC_j Q)(Q^{-1} C_1^{-1}C_2 Q)
= Q^T C_jC_1^{-1}C_2 Q.
\end{equation}

When $j=1,$ by substituting \eqref{Lemma2-0} into \eqref{Lemma2-2}, we have
\begin{equation}\label{Lemma2-3}
(Q^TC_1 Q)J=(Q^TC_1 Q)\cdot\texttt{diag}(\beta_1I_{m_1}, \ldots, \beta_r I_{m_r})=Q^TC_2 Q.
\end{equation}
Since $Q^TC_1Q , Q^TC_2Q $ are both real symmetric and $J$ is a Jordan matrix, Lemma \ref{bd2} asserts that $ Q^TC_1Q $ is a block diagonal matrix with the same partition as $J.$ That is, we can write
\begin{equation}\label{Lemma2-1}
Q^TC_1Q =\texttt{diag}((A_1)_{m_1}, (A_2)_{m_2}, \ldots, (A_r)_{m_r}),
\end{equation}
which proves \eqref{ct1}.
Plugging both \eqref{Lemma2-1} and \eqref{Lemma2-0} into \eqref{Lemma2-3}, we obtain
\begin{eqnarray*}
  &&\texttt{diag}((A_1)_{m_1}, (A_2)_{m_2}, \ldots, (A_r)_{m_r}) \texttt{diag}(\beta_1I_{m_1}, \ldots, \beta_r I_{m_r})\\
  &=& \texttt{diag}(\beta_1A_1, \ldots, \beta_r A_r)=Q^TC_2Q,
\end{eqnarray*}
which proves \eqref{ct1.1}.

Finally, for $j=3,4,\ldots,m$ in \eqref{Lemma2-2}, due to the assumption that
$C_jC_1^{-1}C_2$
are symmetric, so are $Q^T C_jC_1^{-1}C_2 Q$. By Lemma \ref{bd2} again,
 $Q^TC_jQ$ are all block diagonal matrices with the same partition as $J,$ which is exactly \eqref{ct1.2}.
\end{proof}

\begin{rem}\label{rm1}
When there is a non-singular $Q$ that puts $Q^TC_1Q$ and
$Q^TC_2Q$ to \eqref{ct1} and \eqref{ct1.1}, we say that $Q^TC_2Q$ is a non-homogeneous dilation of $Q^TC_1Q$ with scaling factors $\{\beta_1, \beta_2,\ldots,\beta_{r}\}$.
In this case, it is easy to see that $Q^TC_1Q$ and $Q^TC_2Q$ are SDC, say by a congruence $H$. Then, $C_1$ and $C_2$ are
SDC by the congruence $QH$.
\end{rem}	

For $m=2,$ Remark \ref{rm1} and (N1) together give the following result.

\begin{cor}[Greub \cite{Greub}]\label{L0}
Two real symmetric matrices $C_1, C_2,$ with $C_1$ nonsingular,
are SDC if and only if
$C_1^{-1}C_2$ is real similarly diagonalizable.
\end{cor}

It then comes with our first main result, Theorem \ref{dl80}, below.

\begin{thm}\label{dl80}
Let $\mathcal{C}_{ns}=\{C_1, C_2, \ldots, C_m\}\subset \mathcal S^n,~m\ge3$ be a non-singular collection with $C_1$ invertible. Suppose
 for each $i$ the matrix $C_1^{-1}C_i$ is real similarly diagonalizable.
  If $C_jC_1^{-1}C_i$ are symmetric  for $2\leq i<j \leq m,$
  then there always exists a nonsingular real matrix $R$ such that
\begin{align}\label{qt}
R^TC_1R=&\texttt{diag}(A_1, A_2, \ldots, A_s), \nonumber\\
 R^TC_2R=&\texttt{diag}(\alpha^2_1A_1, \alpha^2_2 A_2, \ldots, \alpha^2_ sA_s),\\
  \ldots & \quad \ldots \nonumber\\
  R^TC_m R=&\texttt{diag}(\alpha^m_1A_1, \alpha^m_2 A_2, \ldots, \alpha^m_sA_s),\nonumber
  \end{align}
where $A_t's$ are nonsingular and symmetric, $\alpha^i_t, t=1,2, \ldots, s,$ are real numbers. When the non-singular collection $\mathcal{C}_{ns}$ is transformed into the form of \eqref{qt} by a congruence $R$, the collection $\mathcal{C}_{ns}$ is indeed SDC.
\end{thm}

\begin{proof}
Suppose $C_1^{-1}C_2$ is diagonalized by a non-singular $Q^{(1)}$ with distinct eigenvalues
$\beta^{(1)}_1, \beta^{(1)}_2,\ldots,\beta^{(1)}_{r^{(1)}}$ having multiplicity
$m_1^{(1)},m_2^{(1)},\ldots,m_{r^{(1)}}^{(1)},$ respectively.
Here the superscript $\hbox{}^{(1)}$ denotes the first iteration.
Since $C_jC_1^{-1}C_2$ is symmetric for $j=3,4,\ldots,m,$ Lemma \ref{bd27} assures that
\begin{eqnarray}
C_1^{(1)}={Q^{(1)}}^TC_1{Q^{(1)}}&=& \texttt{diag}\underbrace{((A^{(1)}_1)_{m_1^{(1)}}, (A^{(1)}_2)_{m_2^{(1)}}, \ldots, (A^{(1)}_{r^{(1)}})_{m_{r^{(1)}}^{(1)}})}_{sym.~ \& ~invert.}, \label{ct1-}\\
C_2^{(1)}={Q^{(1)}}^TC_2{Q^{(1)}}&=&\texttt{diag} (\beta^{(1)}_1A^{(1)}_1, \beta^{(1)}_2A^{(1)}_2, \ldots,\beta^{(1)}_{r^{(1)}}A^{(1)}_{r^{(1)}}), \label{ct1-.}\\ %\label{ct1.1-}\\
C_j^{(1)}={Q^{(1)}}^TC_j{Q^{(1)}}&=&\texttt{diag} \underbrace{(C^{(1)}_{j1}, C^{(1)}_{j2},\ldots, C^{(1)}_{jr^{(1)}})}_{sym.},  \enskip j=3,4,\ldots,m; \label{ct1-./} %\label{ct1.2-}
\end{eqnarray}
where all members in $\{C_1^{(1)},C_2^{(1)},C_3^{(1)},\ldots,C_m^{(1)}\}$ adopt the same block structure, each having $r^{(1)}$ diagonal blocks.

As for the second iteration, we use the assumption that $C_1^{-1}C_3$ is similarly diagonalizable. By (ii) of Lemma \ref{L3.0}, \eqref{ct1-}, and \eqref{ct1-./} (for $j=3$), we find that
\begin{equation}\label{step-2.1}
{C_1^{(1)}}^{-1}C^{(1)}_3 = \texttt{diag}\left({A_1^{(1)}}^{-1}C^{(1)}_{31}, \ldots, {A^{(1)}_{r^{(1)}}}^{-1}C^{(1)}_{3r^{(1)}}\right)
\end{equation}
is also similarly diagonalizable. Since a block diagonal matrix is diagonalizable if and only if each of its blocks is diagonalizable, \eqref{step-2.1} implies that each
${A_t^{(1)}}^{-1}C^{(1)}_{3t},~t=1,2,\ldots,r^{(1)}$ is diagonalizable.
Let $Q^{(2)}_t$ (the superscript $\hbox{}^{(2)}$ denotes the second iteration)
diagonalize
${A_t^{(1)}}^{-1}C^{(1)}_{3t}$ into $l_t$ distinct eigenvalues
$\beta^{(2)}_{t1}, \beta^{(2)}_{t2},\ldots,\beta^{(2)}_{tl_t},$ each having multiplicity
$m_{t1}^{(2)},m_{t2}^{(2)},\ldots,m^{(2)}_{tl_t},$ respectively. Then,
$$Q^{(2)}=\texttt{diag}(Q^{(2)}_1, Q^{(2)}_2, \ldots, Q^{(2)}_{r^{(1)}})$$
diagonalizes ${C_1^{(1)}}^{-1}C^{(1)}_3.$

Now, applying Lemma \ref{bd27}
to $\{A_t^{(1)},C^{(1)}_{3t}\}$ for each $t=1,2,\ldots,r^{(1)},$ we have
\begin{eqnarray}
{Q^{(2)}_t}^T A_t^{(1)}{Q^{(2)}_t} &=& \texttt{diag}\underbrace{((A_{t1}^{(2)})_{m_{t1}^{(2)}}, (A_{t2}^{(2)})_{m_{t2}^{(2)}}, \ldots, (A_{tl_t}^{(2)})_{m^{(2)}_{tl_t}})}_{sym.~ \& ~invert.}; \label{ct2-}\\
{Q^{(2)}_t}^T C_{3t}^{(1)}{Q^{(2)}_t}&=& \texttt{diag}(\beta^{(2)}_{t1}A_{t1}^{(2)},\beta^{(2)}_{t2} A_{t2}^{(2)}, \ldots, \beta^{(2)}_{tl_t} A_{tl_t}^{(2)}).\label{ct2-.}
\end{eqnarray}
Let us re-enumerate the indices of all sub-blocks into a sequence from $r^{(1)}$ to $r^{(2)}$:
\begin{eqnarray}
&&\{11,12,\ldots,1l_1\}; \{21,22,\ldots,2l_2\};\cdots;\{r^{(1)}1,r^{(1)}2\ldots,r^{(1)}l_{r^{(1)}}\} \nonumber\\ &\Longrightarrow&\{1,2,\ldots,l_1;l_1+1,l_1+2,\ldots,l_1+l_2; \ldots; \sum_{k=1}^{r^{(1)}-1}l_k+1,\ldots,r^{(2)}\}\label{re-enumerate.1}
\end{eqnarray}
so that
$$A_{11}^{(2)}\rightarrow A_{1}^{(2)};~A_{12}^{(2)}\rightarrow A_{2}^{(2)};\cdots;
A_{1l_1}^{(2)}\rightarrow A_{l_1}^{(2)};~A_{21}^{(2)}\rightarrow A_{l_1+1}^{(2)};~A_{22}^{(2)}\rightarrow A_{l_1+2}^{(2)}\hbox{ and so on.}$$
Assemble \eqref{ct2-} and \eqref{ct2-.} for all $t=1,2,\ldots,r^{(1)}$ together and then use the re-index \eqref{re-enumerate.1}, there is
\begin{eqnarray}
C_1^{(2)}={Q^{(2)}}^TC_1^{(1)}{Q^{(2)}}&=& \texttt{diag}(A^{(2)}_1, A^{(2)}_2, \ldots, A^{(2)}_{r^{(2)}}), \label{ct2}\\
C_3^{(2)}={Q^{(2)}}^TC_3^{(1)}{Q^{(2)}}&=&\texttt{diag} (\beta^{(2)}_1A^{(2)}_1, \beta^{(2)}_2A^{(2)}_2, \ldots,\beta^{(2)}_{r^{(2)}}A^{(2)}_{r^{(2)}}). \label{ct2.2}
\end{eqnarray}
In other words, at the first iteration, $C_1$ is congruent (via $Q^{(1)}$) to a block diagonal matrix $C^{(1)}_1$ of $r^{(1)}$ blocks as in \eqref{ct1-}, while at the second iteration, each of the $r^{(1)}$ blocks is further decomposed (via $Q^{(2)}$) into many more finer blocks ($r^{(2)}$ blocks) as in \eqref{ct2}. Simultaneously, the same congruence matric $Q^{(1)}Q^{(2)}$ makes $C_3$
into $C_3^{(2)}$ in \eqref{ct2.2}, which is a non-homogeneous dilation of $C^{(2)}_1$ with scaling factors $\{\beta^{(2)}_1, \beta^{(2)}_2,\ldots,\beta^{(2)}_{r^{(2)}}\}$.

As for $C_2^{(1)}$ in \eqref{ct1-.}, after the first iteration it has already become a non-homogeneous dilation of $C_1^{(1)}$ in \eqref{ct1-} with scaling factors $\{\beta^{(1)}_1, \beta^{(1)}_2,\ldots,\beta^{(1)}_{r^{(1)}}\}$. Since $C_1^{(1)}$
continues to split into finer sub-blocks as in \eqref{ct2},
$C_2^{(1)}$ will be synchronously decomposed, along with $C_1^{(1)},$ into a block diagonal matrix of $r^{(2)}$ blocks having the original scaling factors $\{\beta^{(1)}_1, \beta^{(1)}_2,\ldots,\beta^{(1)}_{r^{(1)}}\}$. Specifically,
we can expand the scaling factors $\{\beta^{(1)}_1, \beta^{(1)}_2,\ldots,\beta^{(1)}_{r^{(1)}}\}$ to become a sequence of $r^{(2)}$ terms as follows:
\begin{eqnarray}
&&\{\underbrace{\beta^{(1)}_1,\beta^{(1)}_1,\ldots,\beta^{(1)}_1}_{l_1};
\underbrace{\beta^{(1)}_2,\beta^{(1)}_2,\ldots,\beta^{(1)}_2}_{l_2};\cdots;
\underbrace{\beta^{(1)}_{r^{(1)}},\beta^{(1)}_{r^{(1)}},\ldots,
\beta^{(1)}_{r^{(1)}}}_{l_{r^{(1)}}}\}\label{re-enumerate.2}\\
&\triangleq&\{[\beta^{(1)}_1],[\beta^{(1)}_2],\ldots,[\beta^{(1)}_{l_1}];
[\beta^{(1)}_{l_1+1}],\ldots,[\beta^{(1)}_{l_1+l_2}];\ldots;
[\beta^{(1)}_{\sum_{k=1}^{r^{(1)}-1}l_k+1}],\ldots,
[\beta^{(1)}_{r^{(2)}}]\}.\nonumber
\end{eqnarray}
With this notation, we can express
\begin{equation}
C_2^{(2)}={Q^{(2)}}^TC_2^{(1)}{Q^{(2)}}=\texttt{diag} ([\beta^{(1)}_1]A^{(2)}_1, [\beta^{(1)}_2]A^{(2)}_2, \ldots,[\beta^{(1)}_{r^{(2)}}]A^{(1)}_{r^{(2)}}). \label{ct2.1}
\end{equation}

For $C_4^{(1)}$ up to $C_m^{(1)},$ let us take $C_4^{(1)}$ for example because
all the others $C_5^{(1)},C_6^{(1)},$ $\ldots,C_m^{(1)}$ can be analogously taken care of. By the assumption that
$C_4C_1^{-1}C_3$ is symmetric and by (iv) of Lemma \ref{L3.0}, we know that
\begin{equation}\label{step-2.2}
C^{(1)}_4{C_1^{(1)}}^{-1}C^{(1)}_3= \texttt{diag}\left( C^{(1)}_{41} {A_1^{(1)}}^{-1}C^{(1)}_{31},\ldots,C^{(1)}_{4r^{(1)}}  {A_{r^{(1)}}^{(1)}}^{-1}C^{(1)}_{3r^{(1)}}\right)
\end{equation}
is symmetric. Since, for each $t=1,2,\ldots,r^{(1)},$ ${A_t^{(1)}}^{-1}C^{(1)}_{3t}$ is similarly diagonalizable by $Q^{(2)}_t$; and $C^{(1)}_{4t} {A_t^{(1)}}^{-1}C^{(1)}_{3t}$ is symmetric, by Lemma \ref{bd27},
$C_{4t}^{(1)}$ can be further decomposed
into finer blocks to become
\begin{equation}\label{tn}
{Q^{(2)}_t}^T C^{(1)}_{4t}{Q^{(2)}_t}= \texttt{diag}\underbrace{( C_{4,t1}^{(2)} ,  C_{4,t2}^{(2)}, \ldots,   C_{4,tl_t}^{(2)} )}_{sym.}.
\end{equation}
Under the re-indexing formula \eqref{re-enumerate.1} and \eqref{re-enumerate.2}, we have
\begin{equation}\label{tn1}
C_4^{(2)}={Q^{(2)}}^TC_4^{(1)}{Q^{(2)}}=\texttt{diag} (C^{(2)}_{41}, C^{(2)}_{42},\ldots, C^{(2)}_{4r^{(2)}}).
\end{equation}

As the process continues, at the third iteration we use the condition that $C_1^{-1}C_4$ is diagonalizable and
$C_jC_1^{-1}C_4,$ $5\le j\le m$ symmetric to ensure the existence of a congruence $Q^{(3)}$, which puts $\{C_2^{(2)},C_3^{(2)},C_4^{(2)}\}$ as non-homogeneous dilation of the first matrix $C^{(2)}_1,$ whereas from $C_{5}^{(2)}$ up to the last $C_m^{(2)}$ are all block diagonal matrices with
the same pattern as the first matrix $C_1^{(2)}.$  At the final iteration, there is a congruence matrix $Q^{(m-1)}$ that puts
$\{C_2^{(m-1)},C_3^{(m-1)},\ldots,C_m^{(m-1)}\}$
as non-homogeneous dilation of $C^{(m-1)}_1.$ Define
$$R=Q^{(1)}Q^{(2)}Q^{(3)}\cdots Q^{(m-1)}.$$
Then, the non-singular congruence matrix $R$ transforms the collection $\{R^TC_iR:i=1,2,\ldots,m\}$ into block diagonal forms of \eqref{qt}. By Remark \ref{rm1},
the collection $\mathcal{C}_{ns}=\{C_1, C_2, \ldots, C_m\},~m\ge3$ is SDC and the proof is complete.
\end{proof}

With Theorem \ref{dl80}, we can now completely characterize the SDC of a non-singular collection $\mathcal{C}_{ns}=\{C_1, C_2, \ldots, C_m\}.$

\begin{thm}\label{thm1}
Let $\mathcal{C}_{ns}=\{C_1, C_2, \ldots, C_m\}\subset \mathcal S^n,~m\ge3$ be a non-singular collection with $C_1$ invertible.
The collection $\mathcal{C}_{ns}$ is SDC if and only if for each $2\le i \le m,$ the matrix $C_1^{-1}C_i$ is real similarly diagonalizable and $C_jC_1^{-1}C_i, 2\leq i<j \leq m$ are all symmetric.
\end{thm}
	
\section{The SDC problem of singular collection}\label{sec02}

Let $\mathcal{C}_{s}=\{C_1, C_2, \ldots, C_m\}\subset \mathcal S^n$
be a singular collection in which every $C_i\not=0$ is singular.
Consider the first two matrices $C_1,C_2.$ By \cite[Theorem 6]{Jiang}, there is a nonsingular $U_1$ that converts $C_1,C_2$ to block diagonal matrices
{\small \begin{equation}\label{Intro.1}
{\bar C}_1=\texttt{diag}(\underbrace{(C_{11})_{p}}_{invert.~\&~ diag.},0_{n-p});~{\bar C}_2=\texttt{diag}((C_{21})_{p},\underbrace{(C_{26})_{s_1}}_{invert.~\&~ diag.},0_{n-p-s_1})
\end{equation}}where $C_{11}$ and $C_{26}$ are both nonsingular diagonal, $p>0,$ $s_1\ge0;$ and $0_{n-p}$ denotes the zero matrix of
 size $(n-p)\times(n-p).$ Also by \cite[Theorem 6]{Jiang}, the SDC of $\{C_1,C_2\}$ implies the SDC of
$\{(C_{11})_{p},(C_{21})_{p}\}$, the latter of which is a non-singular collection of smaller matrix size $p<n.$
Notice that, in Jiang and Li \cite[Theorem 6]{Jiang}, only the case $s_1>0$ was considered. In this paper, we include the case $s_1=0$ to complete their results. We emphasize that the inclusion of $s_1=0$ complicates the proof by a lot.

Suppose $\{C_{11}, C_{21}\}$ are SDC, say, by $(W)_p$. Let
$Q_1=\texttt{diag}((W)_p, I_{n-p}),$ where $I_{n-p}$ is the identity matrix of dimension $n-p.$ Then, %there is
\begin{eqnarray*}
 {\bar{C}}'_1= Q_1^T\bar{C}_1Q_1&=&\texttt{diag}(\underbrace{ \underbrace{(W^T{C}_{11}W)_{p}}_{\triangleq{\bar{C}}'_{11}:~\small{invert.~\&~ diag.}}, \underbrace{0_{s_1}}_{s_1\ge0}}_{\triangleq{\hat C}_{11}}, 0_{n-p-s_1});
   \\
 {\bar{C}}'_2= Q_1^T\bar{C}_2Q_1&=&\texttt{diag}( \underbrace{(W^T{C}_{21}W)_{p}}_{\triangleq{\bar{C}}'_{21}:\,\small{diag.}}, \underbrace{({C}_{26})_{s_1}}_{\small{\triangleq{\bar{C}}'_{26}:invert.\,\&\, diag.}}, 0_{n-p-s_1}).
\end{eqnarray*}
It allows us to choose a large enough $\mu_1$ such that $\mu_1 {\bar{C}}'_{11}+{\bar{C}}'_{21}$ is invertible (where ${\bar{C}}'_{21}=W^T{C}_{21}W).$ Then,
\begin{eqnarray*}
 \mu_1{\bar{C}}'_1+{\bar{C}}'_2&=& Q_1^T(\mu_1\bar{C}_1+\bar{C}_2)Q_1\nonumber\\
 &=&\texttt{diag}(\underbrace{\underbrace{(\mu_1{\bar{C}}'_{11}
 +{\bar{C}}'_{21})_{p}}_{\small{invert.~\&~ diag.}},\underbrace{({\bar{C}}'_{26})_{s_1}}_{\small{invert.~\&~ diag.}}}_{\triangleq{\hat C}_{21}: \small{invert.~\&~ diag.}}, 0_{n-p-s_1}).
\end{eqnarray*}
Now include $C_3$ for determining the SDC of $\{C_1,C_2,C_3\}$. We first transform $C_3$ by $U_1,$ followed by $Q_1,$ to obtain ${\bar{C}}'_3=Q_1^T(U^T_1C_3U_1)Q_1.$ The idea is to apply \cite[Theorem 6]{Jiang} again to convert
$\mu_1{\bar{C}}'_1+{\bar{C}}'_2$ and ${\bar{C}}'_3$ into the form \eqref{Intro.1}, where, with the help of a sufficiently large $\mu_1>0$, the subblock $({\hat C}_{21})_{p+s_1}$ in $\mu_1{\bar{C}}'_1+{\bar{C}}'_2$ is non-singular and diagonal and thus can be used to determine the SDC of $\{\mu_1{\bar{C}}'_1+{\bar{C}}'_2,{\bar{C}}'_3\}.$
The entire Section 3 is devoted to proving that the idea does indeed work. The main result, Theorem \ref{main}, states that, suppose that the first $n-1$ matrices are SDC (otherwise, it is end of the story), there always exist a sequence of congruences matrices and a sequence of large enough constants which can reduce the SDC of the {\em entire} singular collection $\mathcal{C}_s$ to become the SDC of another non-singular collection $\mathcal{C}_{ns}$ having a smaller matrix size. Then, the result in Section 2 applies.

Let us first modify Lemma 4 in \cite[Jiang and Li]{Jiang}.

\begin{lem}[modified from \cite{Jiang}]\label{lem01}
Let both $C_1,C_2\in \mathcal S^n$ be non-zero singular
with $\text{rank}(C_1)=p<n.$ There exists a nonsingular matrix $U_1,$ which diagonalizes $C_1$ and rearrange its non-zero eigenvalues as
\begin{align}\label{1.}
		\bar{C}_1=U_1^TC_1U_1
		=\left(
		\begin{matrix}
			\underbrace{(C_{11})_p}_{\small{invert.~\&~ diag.}}&0\\ 0&0_{n-p} \end{matrix}
		\right),
		\end{align}
while the same congruence $U_1$ puts $\bar{C}_2=U_1^TC_2U_1$ into two possible forms: either
\begin{align}\label{2.}
		\bar{C}_2=U_1^TC_2U_1
		=\left(
		\begin{matrix}
			(C_{21})_p& C_{22}\\ C^T_{22}&0_{n-p}
		\end{matrix}
		\right),
\end{align}
or
\begin{align}\label{2}
		\bar{C}_2=U_1^TC_2U_1
		=\left(
		\begin{matrix}
			(C_{21})_p&0& C_{25}\\ 0&\underbrace{(C_{26})_{s_1}}_{\small{invert.~\&~ diag.}}&0\\
		C_{25}^T&0& 0_{n-p-s_1}
		\end{matrix}
		\right).
		\end{align}
\end{lem}

Such a $U_1$ in Lemma \ref{lem01} can be obtained by Algorithm 1 \cite[Jiang and Li]{Jiang}. One first finds an orthogonal matrix $Q_1$ such that
\begin{eqnarray}
  \bar{C}_1=Q_1^T C_1Q_1 &=& \texttt{diag}(\underbrace{\texttt{diag}(\alpha_1, \alpha_2, \ldots, \alpha_{p})}_{=(C_{11})_{p},~\small{invert.}},0_{n-p}); \label{to}\\
  Q_1^T C_2Q_1 &=& \left(\begin{matrix}(M_{21})_p&M_{22}\\ M_{22}^T&\underbrace{(M_{23})_{ n-p}}_{sym.}\end{matrix}\right). \label{to1}
\end{eqnarray}
We see that \eqref{to} is already in the form of \eqref{1.}. If $M_{23}=0$ in \eqref{to1},
\begin{equation*}
\bar{C}_2=Q_1^T C_2Q_1=\left(\begin{matrix}(M_{21})_ p&M_{22}\\ (M_{22})^T&0_{n-p}\end{matrix}\right),
\end{equation*}
which is \eqref{2.}.

Otherwise, $\text{rank} M_{23}:=s_1\ge1.$
Let $P_1$ be an orthogonal matrix to diagonalize the symmetric $M_{23}$
as
\begin{equation*}
P_1^TM_{23}P_1=\texttt{diag}(\underbrace{(C_{26})_{ s_1}}_{\small{invert.~\&~ diag.}},0_{n-p-s_1}).
\end{equation*}
Define $H_1=\texttt{diag}(I_p, (P_1)_{n-p})$ and compute
\begin{equation*}
H_1^TQ_1^T C_2Q_1H_1
=\left(
\begin{matrix}	(M_{21})_p&C_{24}&C_{25}\\C_{24}^T&(C_{26})_{s_1}&0\\C_{25}^T&0&0_{n-p-s_1}
			\end{matrix}
		\right),
\end{equation*}
where $(C_{24},C_{25})_{p\times(n-p)}=M_{22}P_1$. Define further that
\begin{equation*}
V_1=\left(\begin{matrix}I_{p}&0&0\\-C_{26}^{-1}C_{24}^T&I_{s_1}&0\\0&0&I_{n-p-s_1}
			\end{matrix}\right);\hbox{ and } U_1=Q_1H_1V_1.
\end{equation*}
Note that the matrix $H_1V_1$ does not change $Q_1^T C_1Q_1$ that we have
\begin{eqnarray*}
\bar{C}_1=U_1^TC_1U_1&=&V_1^TH_1^TQ_1^T C_1Q_1H_1V_1=Q_1^TC_1Q_1~~(\text{as in }\eqref{to}) \\
\bar{C}_2=U_1^TC_2U_1&=&V_1^TH_1^TQ_1^T C_2Q_1H_1V_1 \nonumber\\
&=&\left(
\begin{matrix}
	\underbrace{M_{21}-C_{24}C_{26}^{-1}(C_{24})^T}_{=(C_{21})_p}&0&C_{25}\\0&(C_{26})_{ s_1}&0\\C_{25}^T&0&0_{n-p-s_1}
\end{matrix}
\right).
\end{eqnarray*}
These are what we need in \eqref{1.} and \eqref{2} of Lemma \ref{lem01}.

Jiang and Li \cite{Jiang} assert that:

$\bullet$ ~Suppose $U_1=Q_1H_1V_1$ puts $\bar{C}_1=U_1^TC_1U_1$ and $\bar{C}_2=U_1^TC_2U_1$
into \eqref{1.} and \eqref{2}. Then SDC of $\bar{C}_1$ and $\bar{C}_2$ is equivalent to SDC of the sub-matrices $C_{11}$ (in \eqref{1.}) and $C_{21}$ (in \eqref{2}); and the submatrix $C_{25}$ (in \eqref{2}) is a zero matrix or does not exist.

Here, we would like to add an additional result to supplement Lemma 4 in Jiang and Li \cite{Jiang}:

$\bullet$ Suppose $U_1=Q_1$ puts $\bar{C}_1$ and $\bar{C}_2$ into \eqref{1.} and \eqref{2.}. Then $\bar{C}_1$ and $\bar{C}_2$ are SDC if and only if $C_{11}$ (in \eqref{1.}) and $C_{21}$ (in \eqref{2.}) are SDC; and $C_{22}=0$ (in \eqref{2.}).
The new result will be accomplished by a couple of lemmas below.

\begin{lem}\label{L2-}
Suppose that $A,B\in \mathcal S^n$ of the following forms are SDC %by a congruence $P$ %$A, B$ taking the following formats
\begin{equation}\label{CdA}
A={\rm diag}(\underbrace{(A_1)_p}_{invert.},0_{n-p}),~~
		B		=\left(
		\begin{matrix}
			(B_1)_p&(B_2)_{p\times(n-p)}\\ B_2^T&0_{n-p}
    \end{matrix}
		\right)
\end{equation}
with $A_1$ non-singular and $p<n$. Then, the congruence $P$ can be chosen to be
\begin{equation*}%\label{Dunhill-clip}
P=\left(\begin{matrix}\underbrace{(P_1)_p}_{invert.}&0\\P_3&P_4\end{matrix}\right)~\hbox{ such that }~P^TAP=\left(\begin{matrix}
\underbrace{(P_1^TA_1P_1)_p}_{invert.\&diag.}&0\\0&0_{n-p}\end{matrix}\right)
\end{equation*}
and
\begin{equation*}%\label{Dunhill-clip.}
P^TBP=\left(\begin{matrix}\underbrace{P_1^TB_1P_1+P_1^TB_2P_3+P_3^TB^T_2P_1}_{diag.}&P_1^TB_2P_4
\\ \underbrace{P_4^TB_2^TP_1}_{=0}&0_{n-p}
 \end{matrix}\right)
\end{equation*}
and thus $B$ must be singular. In other words, if $A,B$ take the form
\eqref{CdA} and $B$ is non-singular, then $\{A,B\}$ cannot be SDC.
\end{lem}
\begin{proof} Since $A,B$ are SDC and $\text{rank}(A)=p$ by the assumption, we can choose the congruence $P$ so that the $p$ non-zero diagonal elements of $P^TAP$ are arranged to the north-western corner, while $P^TBP$ is still diagonal. That is,

\begin{equation*}
P=\left(\begin{matrix}(P_{1})_p&P_{2}\\P_{3}&(P_4)_{n-p}\end{matrix}\right)\Longrightarrow
P^TAP=\left(\begin{matrix}\underbrace{(P_1^TA_1P_1)_p}_{invert.~\& ~ diag.}&\underbrace{(P_1^TA_1P_2)_{p\times (n-p)}}_{=\ 0}\\\underbrace{P_2^TA_1P_1}_{=0}&\underbrace{(P_2^TA_1P_2)_{n-p}}_{=\ 0}\end{matrix}\right).
\end{equation*}
Since $P_{1}^TA_1P_1$ is non-singular diagonal and $A_1$ is non-singular, $P_1$ must be invertible. Then, the off-diagonal $P_{1}^TA_1P_2=0$ implies that $P_2=0_{p\times(n-p)}.$ So,
$$P=\left(\begin{matrix}P_1&0\\P_3&P_4\end{matrix}\right)
 \hbox{ and }P^TBP=\left(\begin{matrix}\underbrace{P_1^TB_1P_1+P_1^TB_2P_3+P_3^TB^T_2P_1}_{diag.}&P_1^TB_2P_4
 \\ \underbrace{P_4^TB_2^TP_1}_{=0}&0_{n-p}
 \end{matrix}\right).$$
Notice that $P^TBP$ is singular, and thus $B$ must be singular, too.
The proof is thus complete.
\end{proof}

\begin{lem}\label{L2}
Let $A,B\in \mathcal S^n$ take the following formats
$$A={\rm diag}((A_1)_p,0_{n-p}),~~
		B		=\left(
		\begin{matrix}
			(B_1)_p&(B_2)_{p\times(n-p)}\\ B_2^T&0_{n-p}
    \end{matrix}
		\right),
$$
with $A_1$ non-singular and $B_2$ of full column rank. Then,
$ {\ker} A\cap {\ker} B=\{0\}.$
\end{lem}

\begin{lem}\label{L1}
Let $A,B\in \mathcal S^n$ with ${\ker} A\cap {\ker} B=\{0\}.$ If
$\alpha A+\beta B$ is singular for all real couples $(\alpha, \beta)\in\Bbb R^2,$ then $A, B$ are not SDC.
\end{lem}
\begin{proof}
Suppose contrarily that $A$ and $B$ were SDC by a congruence $P$ such that
 \begin{align*}
 P^TAP=D_1={\rm diag}(a_1, a_2, \ldots, a_n);~ P^TBP=D_2={\rm diag}(b_1, b_2, \ldots, b_n).
 \end{align*}
Then,
$P^T(\alpha A+\beta B)P={\rm diag}(\alpha a_1+\beta b_1,\alpha a_2+\beta b_2, \ldots,\alpha a_n+\beta b_n).$
By assumption, $\alpha A+\beta B$ is singular for all $(\alpha, \beta)\in\Bbb R^2$
so that at least one of $\alpha a_i+\beta b_i=0, \forall (\alpha, \beta)\in\Bbb R^2.$ Let us say $\alpha a_1+\beta b_1=0, \forall (\alpha, \beta)\in\Bbb R^2.$
It implies that $a_1=b_1=0.$ Let $e_1=(1,0,\ldots,0)^T$ be the first unit vector
and notice that $Pe_1\not=0$ since $P$ is non-singular. Then,
$$P^TAPe_1=D_1e_1=0;~P^TBPe_1=D_2e_1=0 \Longrightarrow 0\not=Pe_1\in{\ker} A\cap {\ker}B,$$
which is a contradiction.
\end{proof}

\begin{lem}\label{L2n}
Let $A,B\in \mathcal S^n$ be both singular taking the following formats
\begin{align*}
	A={\rm diag}(\underbrace{(A_1)_p}_{invert.},0_{n-p});	B		=\left(
		\begin{matrix}
			(B_1)_p&B_2\\ B_2^T&0_{n-p}
    \end{matrix}
		\right),
		\end{align*}
with $A_1$ non-singular and $B_2$ of full column-rank. Then $A, B$ are not SDC.
\end{lem}

\begin{proof}
From Lemma \ref{L2}, we know that ${\ker} A\cap {\ker} B=\{0\}.$ If
$\alpha A+\beta B$ is singular for all $(\alpha,\beta)\in\Bbb R^2,$ Lemma \ref{L1}
asserts that $A$ and $B$ are not SDC.
Otherwise, there is $(\tilde\alpha,\tilde\beta)\in\Bbb R^2$
such that $\tilde\alpha A+\tilde\beta B$ is nonsingular. Surely, $\tilde\alpha\ne0,\tilde\beta\ne0.$ Then,
$$C=\dfrac{\tilde\alpha}{\tilde\beta} A+B=\left( \begin{matrix}(\dfrac{\tilde\alpha}{\tilde\beta} A_1+ B_1)_{p}&B_2\\ B_2^T&0\end{matrix}\right) \hbox{ is non-singular }.$$
By Lemma \ref{L2-}, $A$ and $C$ are not SDC. So, $A$ and $B$ are not SDC, either.
\end{proof}

\begin{lem}\label{L2--}
Let $\bar{A}_1,\bar{A}_2,\ldots,\bar{A}_m\in \mathcal S^n$ be singular and $A_1,A_2,\ldots,A_m\in \mathcal S^p,~p<n$ so that
\begin{equation}\label{CdA.}
\bar{A}_i={\rm diag}((A_i)_p,0_{n-p}).~
\end{equation}
Then $\bar{A}_1, \bar{A}_2,\ldots, \bar{A}_m$ are SDC if and only if   $A_1, A_2, \ldots, A_m$ are SDC.
\end{lem}
\begin{proof}
If $A_1, A_2, \ldots, A_m$ are SDC by a congruence $P_1\in \Bbb R^{p\times p}$ then so are $\bar{A}_1,\bar{A}_2,$ $\ldots,\bar{A}_m$ by a congruence $P={\rm diag}((P_1)_{p},I_{n-p}) \in \Bbb R^{n\times n}.$

 Conversely, suppose $\bar{A}_1, \bar{A}_2,\ldots,\bar{A}_m$ are SDC by a nonsingular matrix
    $$P=\left( \begin{matrix}
(P_1)_p& P_2\\
P_3 & (P_4)_{n-p}
\end{matrix}\right).$$
Thus, for every $i=1,2,\ldots,m,$ the matrix  $$P^T \bar{A}_iP =\left( \begin{matrix}
P_1^TA_iP_1& P_1^TA_iP_2\\
P_2^TA_iP_1 & P_2^TA_iP_2
\end{matrix}\right)$$  is diagonal, implying that $P_1^TA_iP_1, P_2^TA_iP_2$ are diagonal.

 Since $P$ is nonsingular, the first $p$ rows of $P$ are linearly independent, i.e.,
 $${\rm rank}(P_1, P_2)=p.$$
  Let $E$ be a permutation matrix that rearranges the  columns of $P$ such that
  $$PE=U=\left( \begin{matrix}
(U_1)_p& U_2\\
U_3 & (U_4)_{n-p}
\end{matrix}\right) $$   with ${\rm rank}( U_1)=p$ and  the matrix  $$U^T \bar{A}_iU =\left( \begin{matrix}
U_1^TA_iU_1& U_1^TA_iU_2\\
U_2^TA_iU_1 & U_2^TA_iU_2
\end{matrix}\right)$$  is still diagonal.  Therefore, $A_1,A_2,\ldots,A_m$ are SDC by congruence $U_1.$
 The proof is complete.
\end{proof}

\begin{lem}\label{L3}
Let $C_1,C_2\in \mathcal S^n$ be both singular and $U_1$ be non-singular that puts $\bar{C}_1=U_1^TC_1U_1$ and $\bar{C}_2=U_1^TC_2U_1$ into \eqref{1.} and \eqref{2.} in Lemma \ref{lem01}.
If $C_{22}$ is nonzero, $\bar{C}_1$ and $\bar{C}_2$ are not SDC.
 \end{lem}
 \begin{proof}
By Lemma \ref{L2n}, if $C_{22}$ is of full column-rank, $\bar{C}_1$ and $\bar{C}_2$ are not SDC. So we suppose that $C_{22}$ has its column rank $q<n-p$ and set $s=n-p-q>0.$ There is a $(n-p)\times (n-p)$ nonsingular matrix $U$ such that
$C_{22}U=\left(
		\begin{matrix}
			\hat{C}_{22}& 0_{p\times s}
		\end{matrix}
		\right),$
where $\hat{C}_{22}$ is  a $p\times q$ full column-rank matrix. Let
$Q={\rm diag}(I_{p},U).$ Then,

\begin{align*}
\tilde{C}_2= Q^T\bar{C}_2Q&=\left(\begin{matrix} I_p& 0_{p\times(n-p)}\\ 0_{(n-p)\times p}& U^T\end{matrix}\right)\left(\begin{matrix} C_{21}& C_{22}\\ C_{22}^T&0\end{matrix}\right)
\left(\begin{matrix} I_p& 0_{p\times(n-p)}\\ 0_{(n-p)\times p}& U\end{matrix}\right) \\
&=
\left(\begin{matrix} C_{21}&\hat{C}_{22}&0_{p\times s}\\ \hat{C}_{22}^T&0_{ q}&0_{q\times s}\\
0_{s\times p}&0_{s\times q}&0_s\end{matrix}\right);
\end{align*}
and
$$\tilde{C}_1=Q^T\bar{C}_1Q=\left(\begin{matrix} C_{11}&0_{p\times q}&0_{p\times s}\\ 0_{q\times p}&0_q&0_{q\times s}\\
0_{s\times p}&0_{s\times q}&0_s\end{matrix}\right).$$
Observe that, by Lemma \ref{L2n}, the two leading principal submatrices
$$A=\left(\begin{matrix} C_{11}&0_{p\times q}\\ 0_{q\times p}&0_q \end{matrix}\right), B=\left(\begin{matrix} C_{21}&\hat{C}_{22}\\ \hat{C}_{22}^T&0_q \end{matrix}\right)$$
of  $\tilde{C}_1$ and $\tilde{C}_2,$ respectively, are not SDC since $C_{11}$ is non-singular (due to \eqref{1.}) and $\hat{C}_{22}$ is of full column rank. By Lemma \ref{L2--}, $\tilde{C}_1$ and $\tilde{C}_2$ cannot be SDC. Then, $\bar{C}_1$ and $\bar{C}_2$ cannot be SDC, either. The proof is complete.
 \end{proof}

In the following, Lemma \ref{L4} and Lemma \ref{lem2} together give a complete set of necessary and sufficient conditions for SDC of two singular matrices $C_1$ and $C_2.$

\begin{lem}\label{L4}
Let $C_1$ and $C_2$ be two symmetric singular matrices of $n\times n.$ Let $U_1$ be the non-singular matrix that puts $\bar{C}_1=U_1^TC_1U_1$ and $\bar{C}_2=U_1^TC_2U_1$ into the format of \eqref{1.} and \eqref{2.} in Lemma \ref{lem01}.
Then, $\bar{C}_1$ and $\bar{C}_2$ are SDC if and only if $C_{11},$ $C_{21}$ are SDC and $C_{22}=0_{p\times r}.$
 \end{lem}

\begin{lem}[\cite{Jiang}, Theorem 6]\label{lem2} Let $C_1$ and $C_2$ be two symmetric singular matrices of $n\times n.$ Let $U_1$ be the non-singular matrix that puts $\bar{C}_1=U_1^TC_1U_1$ and $\bar{C}_2=U_1^TC_2U_1$ into the format of \eqref{1.} and \eqref{2} in Lemma \ref{lem01}.
Then, $\bar{C}_1$ and $\bar{C}_2$ are SDC if and only if $C^{-1}_{11}C_{21}$ is similarly diagonalizable (which amounts to $C_{11}, C_{21}$ being SDC by Corollary \ref{L0}); and $C_{25}$ is a zero matrix or does not exist.
\end{lem}
Suppose $C_1,C_2$ are SDC and we now include $C_3$ to determine the SDC of $\{C_1,C_2,C_3\}.$
By Lemmas \ref{L4} and \ref{lem2}, there is a $U_1$ that converts $C_1,C_2$ to block diagonal matrices
${\bar C}_1=\texttt{diag}((C_{11})_{p},0_{n-p})$ and
${\bar C}_2=\texttt{diag}((C_{21})_{p},(C_{26})_{s_1},0_{n-p-s_1})$
where $C_{11}$ and $C_{26}$ are both non-singular diagonal, but
$s_1\ge0$ could be $0.$ Moreover, SDC of $C_1,C_2$ implies that
$C_{11}, C_{21}$ are SDC, say, by the congruence $(W)_p$.
Let
$Q_1=\texttt{diag}((W)_p, I_{n-p}).$ Then, %there is
\begin{eqnarray}
 {\bar{C}}'_1= Q_1^T\bar{C}_1Q_1&=&\texttt{diag}(\underbrace{ \underbrace{(W^T{C}_{11}W)_{p}}_{\triangleq{\bar{C}}'_{11}:~\small{invert.~\&~ diag.}}, \underbrace{0_{s_1}}_{s_1\ge0}}_{\triangleq{\hat C}_{11}}, 0_{n-p-s_1}); \label{Montblanc}
   \\
 {\bar{C}}'_2= Q_1^T\bar{C}_2Q_1&=&\texttt{diag}( \underbrace{(W^T{C}_{21}W)_p}_{\triangleq{\bar{C}}'_{21}:\,\small{diag.}}, \underbrace{({C}_{26})_{s_1}}_{\small{\triangleq{\bar{C}}'_{26}:invert.\,\&\, diag.}}, 0_{n-p-s_1}).\label{Montblanc.}
\end{eqnarray}
Synchronically, $C_3$ is first transformed to ${\bar{C}}_3$ by $U_1,$ followed by another transformation by $Q_1$ to become
\begin{align}\label{:b}
 {\bar{C}}'_3=Q_1^T\underbrace{U^T_1C_3U_1}_{{\bar{C}}_3}Q_1&=&\left(\begin{matrix}\underbrace{(M_{31})_{p+s_1}}_{\small{sym.},~s_1\ge0}&M_{32}\\M_{32}^T&\underbrace{(M_{33})_{n-p-s_1}}_{\small{sym.}}\end{matrix}\right)
\end{align}

Note that, in \eqref{Montblanc}, ${\bar{C}}'_{11}=W^T{C}_{11}W$ is invertible due to $C_{11}$ being invertible and ${\rm rank}(C_{11})={\rm rank}({\bar{C}}'_{11}).$ It allows us to choose a large enough $\mu_1$ such that $\mu_1 {\bar{C}}'_{11}+{\bar{C}}'_{21}$ is invertible (where ${\bar{C}}'_{21}=W^T{C}_{21}W).$ Then,
\begin{eqnarray}
 \mu_1{\bar{C}}'_1+{\bar{C}}'_2&=& Q_1^T(\mu_1\bar{C}_1+\bar{C}_2)Q_1\nonumber\\
 &=&\texttt{diag}(\underbrace{\underbrace{(\mu{\bar{C}}'_{11}+{\bar{C}}'_{21})_p}_{\small{invert.~\&~ diag.}},\underbrace{({\bar{C}}'_{26})_{s_1}}_{\small{invert.~\&~ diag.}}}_{\triangleq{\hat C}_{21}: \small{invert.~\&~ diag.}}, 0_{n-p-s_1}).\label{Pelikan500}
\end{eqnarray}

Next, we are going to convert the pair $\mu_1{\bar{C}}'_1+{\bar{C}}'_2$ and ${\bar{C}}'_3$ into the form \eqref{1.} and \eqref{2.}; or the form of \eqref{1.} and \eqref{2} in Lemma \ref{lem01}, respectively. Notice that $\mu_1{\bar{C}}'_1+{\bar{C}}'_2=\texttt{diag}(({\hat C}_{21})_{p+s_1},0_{n-p-s_1})$ is already in the form of \eqref{1.}.

$\bullet$ If, in \eqref{:b}, $M_{33}=0,$ ${\bar{C}}'_3$ is thus in the form of \eqref{2.}. Let us rename
\begin{equation}\label{Duofold}
\hat{C}_1={\bar{C}}'_1~ (\hbox{in } \eqref{Montblanc});~\hat{C}_2=\mu_1{\bar{C}}'_1+{\bar{C}}'_2~ (\hbox{in } \eqref{Pelikan500});~
\hat{C}_3={\bar{C}}'_3=\left(\begin{matrix}M_{31}&M_{32}\\M_{32}^T&0\end{matrix}\right)
\end{equation}
and denote their north-west subblocks as in \eqref{Montblanc} and in \eqref{Pelikan500} %the north-west subblocks
\begin{equation}\label{Duofold-1989}
{\hat C}_{11}=\texttt{diag}(({\bar{C}}'_{11})_{p},0_{s_1});~
{\hat C}_{21}=\texttt{diag}(\underbrace{(\mu{\bar{C}}'_{11}+{\bar{C}}'_{21})_p,
({\bar{C}}'_{26})_{s_1}}_{\small{invert.~\&~ diag.}});~{\hat C}_{31}=(M_{31})_{p+s_1} .
\end{equation}
It is easy to see the following result.

\begin{lem}\label{result2..}
Let $\{\hat{C}_1, \hat{C}_2,\hat{C}_3\}$ be singular matrices of the form
\eqref{Duofold}.
Then, $\{\hat{C}_1, \hat{C}_2,\hat{C}_3\}$ are SDC if and only if the north-western sub-blocks of them,$\{{\hat C}_{11}, {\hat C}_{21}, {\hat C}_{31}\},$
as specified by \eqref{Duofold-1989} are SDC; and $M_{32}=0$.
\end{lem}
\begin{proof}
If $M_{32}=0$ and the northwest sub-blocks $\hat{C}_{11},\hat{C}_{21},\hat{C}_{31}$ in \eqref{Duofold-1989} are SDC by $(L_1)_{p+s_1},$ then the matrix $L=\texttt{diag}(L_1,I_{n-p-s_1})$ simultaneously diagonalizes $\hat{C}_1, \hat{C}_2,\hat{C}_3$ via congruence.

Conversely, suppose $\hat{C}_1, \hat{C}_2,\hat{C}_3$
are SDC. In particular, $\hat{C}_2,\hat{C}_3$ are SDC. Since ${\hat C}_{21}$ is non-singular and diagonal whereas $\hat{C}_3$ is in the form of \eqref{2.}, by Lemma \ref{L4}, $M_{32}$ must be $0.$ It implies that $\hat{C}_1, \hat{C}_2,\hat{C}_3$ have the same block structure. Specifically, $\hat{C}_1=\texttt{diag}(({\hat C}_{11})_{p+s_1},0),$
$\hat{C}_2=\texttt{diag}(({\hat C}_{21})_{p+s_1},0),$ and $\hat{C}_3=\texttt{diag}((M_{31})_{p+s_1},0).$ By Lemma \ref{L2--},
$\{{\hat C}_{11}, {\hat C}_{21}, {\hat C}_{31}\}$ are SDC and the proof is complete.
\end{proof}

$\bullet$ Suppose, in \eqref{:b}, $M_{33}\not=0.$ Let an orthogonal $(P_2)_{n-p-s_1}$ be such that

\begin{equation*}
P_2^TM_{33}P_2=\texttt{diag}(\underbrace{(C_{36})_{ s_2}}_{\small{invert.~\&~ diag.,~s_2>0}},0_{n-p-s_1-s_2}),
\end{equation*}
with which we can form $H_2=\texttt{diag}(I_{p+s_1}, P_2)$ and compute
\begin{equation}\label{c-}
H_2^T{\bar{C}}'_3H_2
=\left(
\begin{matrix}	(M_{31})_{p+s_1}&C_{34}&C_{35}\\C_{34}^T&(C_{36})_{s_2}&0\\C_{35}^T&0&0_{n-p-s_1-s_2}
			\end{matrix}
		\right),
\end{equation}
where $(C_{34},C_{35})_{p\times(n-p)}=M_{32}P_2$. Define further that
\begin{equation}\label{syu.-}
V_2=\left(\begin{matrix}I_{p+s_1}&0&0\\-C_{36}^{-1}C_{34}^T&I_{s_2}&0\\0&0&I_{n-p-s_1-s_2}
			\end{matrix}\right),~\hbox{and}~U_2=H_2V_2
\end{equation}
so that

\begin{align}\label{M400}
\breve{C}_3\triangleq U_2^T{\bar{C}}'_3 U_2=\left(
\begin{matrix}
	\underbrace{M_{31}-C_{34}C_{36}^{-1}C_{34}^T}_{\triangleq({\breve C}_{31})_{p+s_1},~{sym.}}&0&C_{35}\\0&\underbrace{(C_{36})_{s_2}}_{\small{invert.~\&~ diag.}}&0\\C_{35}^T&0&0_{n-p-s_1-s_2}
\end{matrix}
\right).
\end{align}
More importantly, the transformation $U_2,$
\begin{equation}\label{U.2}
U_2=H_2V_2=\left(\begin{matrix}
			I_{p+s_1}&0\\ -P_2\begin{bmatrix}C_{36}^{-1}C_{34}^T\\0  \end{bmatrix}&(P_2)_{n-p-s_1} \end{matrix}
		\right),
\end{equation}
does not change ${\bar{C}}'_1$ in \eqref{Montblanc} and $\mu_1{\bar{C}}'_1+{\bar{C}}'_2$ in \eqref{Pelikan500}, in the sense that
\begin{align}
\breve{C}_1\triangleq U_2^T{\bar{C}}'_1 U_2={\bar{C}}'_1=\texttt{diag}( \underbrace{(W^T{C}_{11}W)_p, 0_{s_1}}_{\triangleq \breve{C}_{11}=\hat{C}_{11}:~\small{diagonal}}, 0_{n-p-s_1}) \label{parker75}
\end{align}
\begin{eqnarray}
 \breve{C}_2&=&U_2^T(\mu_1\bar{C}'_1+\bar{C}'_2)U_2=\mu_1\bar{C}'_1+\bar{C}'_2\nonumber\\
 &=&\texttt{diag}(\underbrace{(\mu_1\bar{C}'_{11}+\bar{C}'_{21})_{p},(\bar{C}'_{26})_{s_1}}_{\triangleq \breve{C}_{21}=\hat{C}_{21}:~\small{invert.~\&~ diag.}}, 0_{n-p-s_1}).\label{Pelikan500.}
\end{eqnarray}
Notice that, in \eqref{parker75} and \eqref{Pelikan500.}, $\hat{C}_{11}$ is renamed as $\breve{C}_{11}$, while $\hat{C}_{21}$ becomes $\breve{C}_{21}.$
Then, we have the following main result.

\begin{lem}\label{result2.}
The singular collection $\{\breve{C}_1, \breve{C}_2,\breve{C}_3\}$ in \eqref{parker75}, \eqref{Pelikan500.}, in \eqref{M400} are SDC if and only if the north-western sub-blocks of them, i.e. $\{{\breve C}_{11}, {\breve C}_{21}, {\breve C}_{31}\},$ are SDC; and $C_{35}$ in \eqref{M400} is a zero matrix or does not exist.

\end{lem}

\begin{proof} The sufficiency of Lemma \ref{result2.} is easy.
If $C_{35}$ in \eqref{M400} is a zero matrix or does not exist, and if the northwest sub-blocks $\breve{C}_{11},\breve{C}_{21},\breve{C}_{31}$ in \eqref{parker75}, \eqref{Pelikan500.} and \eqref{M400}
are SDC by $(L_1)_{p+s_1},$ then the matrix $L=\texttt{diag}(L_1,I_{s_2},I_{n-p-p_2-s_2})$ simultaneously diagonalizes $\{\breve{C}_1, \breve{C}_2,\breve{C}_3\}$ via congruence.

To prove the necessity, suppose that $\breve{C}_1, \breve{C}_2,\breve{C}_3$ are SDC by a congruence matrix $Q.$ In particular, $\breve{C}_2,\breve{C}_3$ are SDC in which $$\breve{C}_{21}=\texttt{diag}(\underbrace{(\mu\bar{C}'_{11}+\bar{C}'_{21})_{p},(\bar{C}'_{26})_{s_1}}_{\small{invert.~\&~ diag.}})~(\hbox{in }\breve{C}_2),~\underbrace{(C_{36})_{s_2}}_{\small{invert.~\&~ diag.}}~(\hbox{in }\breve{C}_3)$$ are non-singular diagonal. By Lemma \ref{lem2},
two matrices (here they are $\breve{C}_2,\breve{C}_3$) in the form of \eqref{1.} and \eqref{2} are SDC, there must be $C_{35}=0$ in $\breve{C}_3$ \eqref{M400} or $C_{35}$ does not exist. Let us assume that $C_{35}=0.$ Then,
\begin{eqnarray*}
\breve{C}_{2}&=&\texttt{diag}(\underbrace{(\breve{C}_{21})_{p+s_1}}_{\small{invert.~\&~ diag.}},0_{s_2},0_{n-p-s_1-s_2}) \label{146G}\\
\breve{C}_{3}&=&\texttt{diag}((\breve{C}_{31})_{p+s_1},\underbrace{({C}_{36})_{s_2}}_{\small{invert.~\&~ diag.}},0_{n-p-s_1-s_2}) \label{146G.}
\end{eqnarray*}
where $\breve{C}_{31}=M_{31}-C_{34}C_{36}^{-1}C_{34}^T$ has been defined in \eqref{M400}.

By Lemma \ref{L2--}, two matrices (which are $\breve{C}_{2}, \breve{C}_{3}$ with $p=p+s_1+s_2,$)
of form \eqref{CdA.} are SDC,
the congruence $Q$ that diagonalizes $\breve{C}_{2},\breve{C}_{3}$ can be chosen to be
\begin{equation}
Q=\left(\begin{matrix}(Q_1)_{p+s_1}&(Q_2)_{(p+s_1)\times s_2}&0\\ (Q_3)_{s_2\times (p+s_1)}&(Q_4)_{s_2}&0\\0&0&I_{n-p-s_1-s_2}
		\end{matrix}\right)\label{Q.}
\end{equation}
such that the first $p+s_1$ diagonal entries of the diagonal matrix $Q^T\breve{C}_{2}Q$ are all non-zero. We shall show that $Q_2=0_{(p+s_1)\times s_2}$ and $Q_3=0_{s_2\times (p+s_1)}$ so that
$Q=\texttt{diag}\left((Q_1)_{p+s_1},(Q_4)_{s_2},I_{n-p-s_1-s_2}\right).$

By $Q$ in \eqref{Q.}, $\breve{C}_{2}$ is congruent to the diagonal matrix
		$$Q^T\breve{C}_{2}Q=\left(\begin{matrix}(Q_1^T \breve{C}_{21}Q_1)_{p+s_1}&Q_1^T \breve{C}_{21}Q_2&0\\
		Q_2^T \breve{C}_{21}Q_1&(Q_2^T \breve{C}_{21}Q_2)_{s_2}&0\\
		0&0&0_{n-p-s_1-s_2}
		\end{matrix}\right)$$
in which $Q_1^T {\breve C}_{21}Q_1$ is nonsingular diagonal. Since $\breve C_{21}$ is also non-singular, it implies that $Q_1$ must be nonsingular. Then, due to the off-diagonal block $Q_1^T \hat{C}_{21}Q_2=0,$ we see that $Q_2=0.$ Then, %Moreover, $Q_4$ must be non-singular, since $Q$ is non-singular. Then,
\begin{equation}
Q^T\breve{C}_{2}Q=\texttt{diag}(\underbrace{(Q_1^T \breve{C}_{21}Q_1)_{p+s_1}}_{\small{invert.~\&~ diag.}},0_{s_2},0_{n-p-s_1-s_2}).\label{146G..}
\end{equation}

Since $\breve{C}_1$ in \eqref{parker75} and $\breve{C}_{2}$ in \eqref{Pelikan500.}
adopt the same block structure, there also is
\begin{equation}
Q^T\breve{C}_{1}Q=\texttt{diag}(\underbrace{(Q_1^T \breve{C}_{11}Q_1)_{p+s_1}}_{\small{diag.}},0_{s_2},0_{n-p-s_1-s_2}).\label{.146G.}
\end{equation}

The same congruence $Q$ also diagonalizes $\hat{C}_3.$ Since $Q_2=0$ in \eqref{Q.},
\begin{align}\label{5}
		  Q^T\breve{C}_3Q
		  =\left(
	  \begin{matrix}
	  	Q_1^T\breve{C}_{31}Q_1+Q_3^TC_{36}Q_3&Q_3^TC_{36}Q_4&0\\ Q_4^TC_{36}Q_3&Q_4^TC_{36}Q_4&0\\
	0&0&0
		\end{matrix}\right) \hbox{ is diagonal }
\end{align}
so that $Q_4^TC_{36}Q_3=0.$ From \eqref{Q.}, since $Q$ is non-singular and we have known that $Q_2=0,$ there must be $Q_4$ non-singular. From \eqref{M400}, we also know that
$C_{36}$ is nonsingular. Then, $Q_4^TC_{36}Q_3=0$ implies that $Q_3=0,$ which proves that the congruence  $Q=\texttt{diag}\left((Q_1)_{p+s_1},(Q_4)_{s_2},I_{n-p-s_1-s_2}\right)$ and \eqref{5} becomes
\begin{equation}%\label{7.1}
Q^T\breve{C}_3Q
		=\texttt{diag}(\underbrace{(Q_1^T \breve{C}_{31}Q_1)_{p+s_1}}_{\small{diag.}},\underbrace{(Q_4^TC_{36}Q_4)_{s_2}}_{\small{diag.}},0_{n-p-s_1-s_2}).\label{149.}
\end{equation}
Combining \eqref{.146G.},\eqref{146G..},\eqref{149.}, we see that
if $\breve{C}_1, \breve{C}_2,\breve{C}_3$ are SDC by $Q,$ then the north-western blocks
$\breve{C}_{11},\breve{C}_{21},\breve{C}_{31}$ are SDC by $Q_1.$ The proof is complete.
\end{proof}

In summary, when $C_1$ and $C_2$ are SDC, there is
$${\breve C}_1=U_2^T\overbrace{Q^T_1\underbrace{U^T_1(C_1)U_1}_{=\bar{C}_1\rm in\ Lemma\ \ref{lem01}}Q_1}^{=\bar{C}'_1\rm\ in\ \eqref{Montblanc}}U_2$$
where $U_1$ is from Lemma \ref{lem01} that puts $C_1,C_2$ in the form of  $\bar{C}_1=U_1^TC_1U_1$  \eqref{1.} and $\bar{C}_2=U_1^TC_2U_1$ \eqref{2}; while
$Q_1$ from \eqref{Montblanc}-\eqref{Montblanc.} diagonalizes simultaneously $\bar{C}_1$ and $\bar{C}_2;$ finally $U_2$ from \eqref{U.2} puts
$\bar{C}'_3=Q^T_1U^T_1C_3U_1Q_1$ in the form of \eqref{M400}.
In addition,
$${\breve C}_2=U_2^TQ^T_1U^T_1(\mu C_1+C_2)U_1Q_1U_2;~{\breve C}_3=U_2^TQ^T_1U^T_1(C_3)U_1Q_1U_2.$$
It is obvious that ${\breve C}_1,{\breve C}_2,{\breve C}_3$ are SDC if and only if $C_1,\mu C_1+C_2,C_3$ are SDC; and, if and only if $C_1,C_2,C_3$ are SDC. Therefore, from Lemma \ref{result2.}, we have the following result.
%Lemma \ref{result2..} and
\begin{lem}\label{result3..}
Let $\{C_1,C_2,C_3\}\subset\mathcal{S}^n$ be a singular collection and assume that $C_1,C_2$ are SDC. Then,
there is a non-singular $U$ and a constant $\mu$ such that
$\breve{C}_1=U^TC_1U, \breve{C}_2=U^T(\mu C_1+C_2)U,\breve{C}_3=U^TC_3U$ be singular matrices of the forms \eqref{parker75}, \eqref{Pelikan500.} and \eqref{M400}, respectively.
Moreover, the collection $\{C_1,C_2,C_3\}$ is SDC if and only if
the northwestern nonsingular subblocks of them, $\{{\breve C}_{11}, {\breve C}_{21}, {\breve C}_{31}\},$ are SDC; and $C_{35}$ in \eqref{M400} is either zero or does not exist.
\end{lem}

Lemmas \ref{result2..} and \ref{result3..} can be easily extended to more than three matrices. Theorem \ref{main} below can be proved by induction, but the detail is messy. So we omit its proof.

\begin{thm}\label{main}
Let $\mathcal{C}_s=\{C_1, C_2,\ldots, C_m\}\subset \mathcal{S}^n,~m\ge3$ be a singular collection in which none is zero. If $C_1, C_2,\ldots, C_{m-1}$ are SDC, then there exist a nonsingular real matrix $Q$ and a positive vector
$\mu=(\mu_1,\mu_2,\ldots,\mu_{m-2}, 1)\in\Bbb R^{m-1}_{++}$ such that %$C_i$ is congruent to $\bar{C}_i$ under $Q$ as follows:
\begin{eqnarray}
  \bar{C}_1 &=& Q^TC_1Q=\texttt{diag}((C_{11})_p,0_{n-p}),~~p<n; \nonumber \\ %\label{c1} \\
  \bar{C}_2 &=& Q^T(\mu_1C_1+C_2)Q=\texttt{diag}((C_{21})_p,0_{n-p}); \nonumber \\%\label{c2} \\
  \bar{C}_3 &=& Q^T(\mu_2(\mu_1C_1+C_2)+C_3)Q=\texttt{diag}((C_{31})_p,0_{n-p}); \nonumber \\ %\label{c3}\\
   &\vdots&  \nonumber\\
  \bar{C}_{m-1} &=& Q^T(\mu_{m-2}( \cdots\mu_3(\mu_2(\mu_1C_1+C_2)+C_3)+C_4)+\cdots +C_{m-2})+ C_{m-1} )Q  \nonumber\\
  &=& \texttt{diag}((C_{(m-1)1})_p,0_{n-p});\label{6}
\end{eqnarray}
  and either
  \begin{eqnarray}
  \bar{C}_m &=& Q^TC_mQ = \left(\begin{matrix}(C_{m1})_p&C_{m2}\\C_{m2}^T&0_{n-p}
			\end{matrix}\right); \label{7}
\end{eqnarray}
or
\begin{eqnarray}
  \bar{C}_m &=& Q^TC_mQ = \left(\begin{matrix}(C_{m1})_p&0&C_{m5}\\0&(C_{m6})_s&0
			\\ C_{m5}^T&0&0_{n-p-s} \end{matrix}\right),~s\le n-p, \label{7.}
\end{eqnarray}
where
\begin{itemize}
\item the sub-matrices $(C_{i1})_p, i=1,2,\ldots,m-1,$ are all diagonal of the same size. In particular, $(C_{(m-1)1})_p$ in \eqref{6} is nonsingular;
\item in \eqref{7}, $(C_{m1})_p$ is symmetric;
\item in \eqref{7.}, $(C_{m1})_p$ is symmetric, $(C_{m6})_s$ is  nonsingular diagonal; $C_{m5}$ is either a $p\times (n-p-s)$ matrix if $s<n-p$ or does not exist if $s=n-p.$
\end{itemize}
Moreover, the following three statements are equivalent.
\begin{itemize}
  \item[(i)] all matrices in the collection $\mathcal{C}_s$ are SDC;
  \item[(ii)] all matrices in the collection $\bar{\mathcal{C}}_s=\{\bar C_1,\bar C_2,\ldots, \bar C_m\}$ are SDC;
  \item[(iii)] either sub-blocks $C_{11}, C_{21}, \ldots, C_{m1}$ with $C_{m1}$ coming from \eqref{7} are SDC and $C_{m2}=0;$ or sub-blocks $C_{11}, C_{21}, \ldots, C_{m1}$ with $C_{m1}$ coming from \eqref{7.} are SDC and either $C_{m5}=0$ or $C_{m5}$ does not exist.
\end{itemize}
\end{thm}

\section{Implications of SDC in Quadratic Optimization}\label{sec3}
	
The QCQP has the following formulation with $m$ the number of constraints:
	\begin{equation*}
	{\rm(P_m)} \hspace*{0.6cm}
	\begin{array}{lll}
	\lambda^*=&\min &f_0(x)=x^TQ_0x+q_0^Tx\\
	&{\rm s.t.} & f_i(x)=x^TQ_ix+q_i^Tx+a_i\le 0,~ i=1,2,\ldots,m,
	\end{array}
	\end{equation*}
where $Q_i\in \mathcal{S}^n$, $q_i\in\Bbb R^n$ and $a_i\in\Bbb R.$

%If $m=1$ and $Q_1\succ0,$ ${\rm(P_1)}$ is then referred as the trust region subproblem (TRS). See a review paper by Yuan \cite{Yuan}.
%If $m=2$ and $Q_1\succ0, Q_2\succeq0,$ ${\rm(P_2)}$ is called the CDT problem posed by Celis, Denis, Tapia \cite{CDT}.
%%
%%Both (TRS) and (CDT) can be solved in polynomial time, despite of the (possibly) non-convex objective function ($Q_0$ may be indefinite). For polynomial solvability of (CDT), the results are rather recent. Please refer to Bienstock \cite{Bienstock16};
%%Sakaue et al. \cite{Sakaue-Nakatsukasa-Takeda-Iwata16}; and Consolini and Locatelli \cite{Consolini-Locatelli}.
%In general, QCQP is NP-hard even when $Q_0$ has just one negative eigenvalue and $Q_1=Q_2= \cdots =Q_m=0.$ See \cite{PP}.

With the iterative procedure developed in Sections 2 and 3, we can now
determine whether or not
$Q_0, Q_1,\ldots, Q_m$ are SDC in polynomial time. If affirmative,
a congruence matrix $P$ is obtained so that
	$$P^TQ_iP={\rm diag}(\alpha^i_1,\ldots,\alpha^i_n).$$
By change of variables $x=Py,$ the quadratic forms $x^TQ_ix$ become the sums of squares in $y.$ That is,
$$x^TQ_ix=y^TP^TQ_iPy=\sum_{j=1}^n\alpha^i_jy_j^2.$$

	If ${\rm(P_m)}$ is homogeneous, i.e., $q_i=0, i=0,1,\ldots,m,$ by $z_j=y_j^2,$ ${\rm(P_m)}$ is reduced to a linear programming in non-negative variables $z:$
	\begin{equation*}
	{\rm(P_m)} \hspace*{0.6cm}
	\begin{array}{lll}
	\lambda^*=&\min &\sum_{j=1}^n\alpha^0_jz_j\\
	&{\rm s.t.} & \sum_{j=1}^n\alpha^i_jz_j+ a_i\le 0,~ i=1,2,\ldots,m,\\
	& & z_j\ge0, ~j=1,2,\ldots,n.
	\end{array}
	\end{equation*}
So, we can see that the biggest beneficiary of being SDC is homogeneous quadratic optimization. The benefit can be extended to cover the sum-of-generalized-Rayleigh-quotients problem (see, for example, \cite{Bong2,Hong,Sun,WangXia})
	\begin{align*}\label{sum-of-two}
	{\rm(R_m) \hspace*{1cm}}  \max_{x\in\mathbb{R}^n\setminus\{0\}}
	\left\{\dfrac{x^TA_1x}{x^TB_1x}+\dfrac{x^TA_2x}{x^TB_2x}+\ldots +\dfrac{x^TA_mx}{x^TB_mx}\right\}
	\end{align*}
where $A_i, B_i\in{S}^n$ and $B_i\succ0.$
When $A_1,A_2,\ldots, A_m;~B_1, B_2, \ldots,B_m$ are SDC, problem
${\rm (R_m)}$ is reduced to maximizing the sum-of-linear-ratios
	$${\rm (SLR_m)}\hspace*{0.8cm} \max_{z\ge 0, z\ne0}\sum_{i=1}^m\dfrac{\alpha_i^Tz}{\beta_i^Tz}.$$
Even though both ${\rm (R_m)}$ and ${\rm (SLR_m)}$ are NP-hard,
the latter can be better approximated by some methods,
 such as an interior algorithm in \cite{Fre}, a range-space approach in \cite{Sheu-Wu-Ilker} and a branch-and-bound algorithm in \cite{Kuno,Jiao}. Please see a good survey on sum-of-ratios problems in \cite{Schaible}.

The difficulty arises when ${\rm(P_m)}$ is nonhomogeneous. Even with the SDC condition, by setting
$\alpha_i=(\alpha^i_1,\ldots,\alpha^i_n)^T,$ $\xi_i=P^Tq_i$ and
$z_j=y_j^2, j=1,2,\ldots,n,$ ${\rm(P_m)}$ is still a non-convex quadratic programming problem:
\begin{equation*}
	{\rm(P_m)} \hspace*{0.6cm}
	\begin{array}{lll}
	\lambda^*=&\min &f_0(y,z)=\alpha_0^Tz+\xi_0^Ty\\
	&{\rm s.t.} & f_i(y,z)=\alpha_i^Tz +\xi_i^Ty+ a_i\le 0,~ i=1,2,\ldots,m,\\
	& & y_j^2= z_j, j=1,2,\ldots,n
	\end{array}
	\end{equation*}
with non-convex equality constraints $y_j^2= z_j.$ Ben-Tal and Hertog \cite{Ben-Tal} showed that, for $m=1$ and a few special cases of $m=2,$ by relaxing
$y_j^2=z_j$ to become $y_j^2\le z_j, j=1,2,\ldots,n,$ the relaxed convex second-order cone programming (SOCP) problem is exact. However, this relaxation skill cannot be extended to more general cases.

The implication of Ben-Tal and Hertog's result in \cite{Ben-Tal} indicates that the SDC property along is not sufficient to resolve the difficult non-homogeneous QCQP problem for $m\ge2.$ Information from the convexity of the joint numerical range
$$
\mathcal{R}(f_0,f_1,\ldots,f_m(x))_{\mathbb{R}^n}=\big\{ \big( f_0(x),f_1(x),\ldots,f_m(x) \big) \in \mathbb{R}^m ~\big|~x \in \mathbb{R}^n \big\}
$$
or
$$
\mathcal{R}(f_0,f_1,\ldots,f_m(x))_{\|x\|=1}=\big\{ \big(f_0(x),f_1(x),\ldots,f_m(x) \big) \in \mathbb{R}^m ~\big|~ \|x\|=1,~x \in \mathbb{R}^n \big\}
$$
is expected to play important roles. Our suspicion is based on the following two results in literature. The first one comes from Theorem 2.1 in \cite{Polyak98} by Polyak, which states that, suppose $Q_0,Q_1,Q_2\in\mathcal{S}^n$ with $n\ge3.$ Then, the following two statements are equivalent:
\begin{itemize}
\item[]$({\rm FE_1})$~$(x^TQ_0x=0) \wedge (x^TQ_1x = 0)\wedge (x^TQ_2x = 0)~\Longrightarrow~ (x=0);$ and \\ the joint range set $\big\{ \big( x^TQ_0x,x^TQ_1x,x^TQ_2x \big) \in \mathbb{R}^3 ~\big|~x \in \mathbb{R}^n \big\}$ is a pointed closed convex cone;
\item[]$({\rm FE_2})$~($\exists \mu_0,\mu_1,\mu_2 \in \mathbb{R}$)~$\mu_0 Q_0 + \mu_1Q_1+\mu_2 Q_2\succ0.$
\end{itemize}
The second result is from Corollary 1.1 in \cite{Jacobson} by Jacobson. It states that, if $Q_0,Q_1,\ldots,Q_r;$ $Q_{r+1},\ldots,Q_m\in\mathcal{S}^n$ are SDC. Then, the following two statements are equivalent:
\begin{itemize}
\item[]$({\rm FE_3})$~$(x^TQ_ix=0,i=1,\ldots,r)\wedge(x^TQ_jx\ge0,j=r+1,\ldots,m)$\\
    $\Longrightarrow~ x^TQ_0x\ge0;$
\item[]$({\rm FE_4})$~($\exists \mu_i,i=1,\ldots,r \in \mathbb{R};~\mu_j\ge0,j=r+1,\ldots,m\in \mathbb{R}$)\\ $Q_0-\sum_{i=1}^{r}\mu_i Q_i -\sum_{j=r_1}^{m} \mu_jQ_j\succeq0.$
\end{itemize}
Both equivalence: $({\rm FE_1})\sim({\rm FE_2})$ and $({\rm FE_3})\sim({\rm FE_4})$ can be viewed as the extension of Finsler-Calabi-type theorem to more than two symmetric matrices.
Both equivalence: $({\rm FE_1})\sim({\rm FE_2})$ and $({\rm FE_3})\sim({\rm FE_4})$ also involve the SDC condition, as, when
$Q_0,Q_1,Q_2$ are SDC, the joint range set $\big\{ \big( x^TQ_0x,x^TQ_1x,x^TQ_2x \big) \in \mathbb{R}^3 ~\big|~x \in \mathbb{R}^n \big\}$ is a closed convex polyhedral cone, though not necessarily pointed. The difference between $({\rm FE_2})$ and $({\rm FE_4})$ is positive-definiteness and positive semi-definiteness, which represents, in optimization, the sufficient condition by the former and the necessary condition by the latter. If, for $m=2$ at least, the gap is not large between $({\rm FE_2})$ and $({\rm FE_4})$, we wish that we will be able to say something about the convexity of the non-homogeneous joint numerical range
$$
\mathcal{R}(f_0,f_1,f_3)_{\mathbb{R}^n}=\big\{ \big( x^TQ_0x+q_0^Tx,x^TQ_1x+q_1^Tx,x^TQ_2x+q_2^Tx \big) \in \mathbb{R}^m ~\big|~x \in \mathbb{R}^n \big\}
$$
under the SDC condition; and also resolve the (non-homogeneous) QCQP for $m=2$ under the SDC condition of the quadratic forms.

\end{document}